\documentclass{article}
\usepackage{authblk}
\usepackage[utf8]{inputenc}
\usepackage{amsmath}
\usepackage{amsbsy}
\usepackage{amssymb}
\usepackage{graphicx}
\usepackage{dsfont}
\usepackage{upgreek}
\usepackage{textcomp}
\usepackage{braket}
\usepackage[margin=1in]{geometry}
\usepackage{amsthm}
\usepackage{mathrsfs}
\usepackage{mathtools}
\usepackage[table,svgnames]{xcolor}
\usepackage{graphicx}
\usepackage{tikz}
\usepackage{float}
\usepackage{enumerate}
\usetikzlibrary{arrows}
\usepackage[toc,page]{appendix}
\usepackage{cleveref}
\usepackage{xcolor}

\crefname{appsec}{appendix}{appendices}
\numberwithin{equation}{section}

\newtheorem{theorem}{Theorem}[section]
\newtheorem{lemma}{Lemma}[section]

\newtheorem{remark}{Remark}[section]
\newtheorem{definition}{Definition}[section]
\newtheorem{proposition}{Proposition}[section]

\providecommand{\keywords}[1]
{
  \textbf{Key words.} #1
}
\providecommand{\classification}[1]
{
  \textbf{Mathematics Subject Classification.} #1
}
\title{Small Time Behavior and Summability for the Schr\"odinger Equation}
\author{Brian Choi}
\affil{\footnotesize{Department of Mathematics and Statistics, 
  Boston University, 111 Cummington Mall, Boston, 02215, MA, USA\\choigh@bu.edu}}
\date{}
\begin{document}
\maketitle
\begin{abstract}
We consider the Carleson's problem regarding small time almost everywhere convergence to initial data for the Schr\"odinger equation, both linear and nonlinear on $\mathbb{R}$. It is shown, via the smoothing effect of the Schr\"odinger flow, that the (sharp) result proved by Dahlberg and Kenig for initial data in Sobolev spaces still holds when one considers the full Schr\"odinger equation with a certain class of potentials. As for $s<\frac{1}{4}$, the failure of $L^p$-boundedness of the (localized) maximal operator is investigated.
\end{abstract}
\\
\classification{40A30,28A20,35Q55,42A24}
\section{Introduction.}
\label{}
Consider the Cauchy problem
\begin{equation*}
      \begin{cases}
    i\partial_t u + \Delta u =0,  (x,t) \in \mathbb{R}^n \times \mathbb{R}\\
    u(0) = u_0\in L^2(\mathbb{R}^n)
    \end{cases}
\end{equation*}
noting that  $-\Delta$ is a non-negative operator. A straightforward computation with the Fourier transform yields
\begin{equation*}
    u(x,t)= e^{it\Delta}u_0 (x) = \int \hat{u}_0(\xi)e^{-it|\xi|^2 + i\xi \cdot x}d\xi,
\end{equation*} 
where
\begin{equation*}
    \hat{f}(\xi) = \int f(x)e^{-ix\cdot\xi}dx,\: \Check{F}(x) = (2\pi)^{-n}\int F(\xi)e^{ix \cdot \xi}d\xi.
\end{equation*}
In this paper we continue to build upon a question initially posed by \cite{carleson1966convergence}: what is the minimal Sobolev regularity $s_\ast$ for which  $e^{it\Delta}f \xrightarrow[t\to 0]{} f$ almost everywhere (a.e.) with respect to the Legesgue measure, for all $f \in H^{s_\ast}(\mathbb{R})$? Carleson originally proved a positive result, that any $f \in H^s(\mathbb{R})$ for $s \geq \frac{1}{4}$ exhibits almost everywhere (a.e.) convergence. Soon \cite{dahlberg1982note} showed that Carleson's result is sharp. In higher dimensions, this problem is closed except at the endpoint $s=\frac{n}{2(n+1)}$. In $n=2$, \cite{du2017sharp} showed sufficiency for $s>\frac{1}{3}$ while \cite{sjolin1987regularity} and \cite{vega1988schrodinger} independently showed sufficiency in $n \geq 3$ for $s>\frac{1}{2}$. \cite{bourgain2013schrodinger} showed sufficiency for $s > \frac{2n-1}{4n}$ for $n \geq 2$, and though it had long been believed that $s>\frac{1}{4}$ is the sharp sufficient condition in higher dimensions, \cite{bourgain2016note} showed necessity for $s \geq \frac{n}{2(n+1)}$ in $n \geq 2$. Recently \cite{du2018pointwise} showed sufficiency for $s>\frac{n+1}{2(n+2)}$ for $n\geq 3$, which was subsequently improved to the sharp condition $s>\frac{n}{2(n+1)}$ by \cite{du2018sharp}. Many of these results generalise nicely to $i\partial_t u + \Phi(D)u=0$ where $\Phi$ is a Fourier multiplier satisfying $|D^\gamma \Phi (\xi)| \lesssim |\xi|^{\alpha-|\gamma|}$ and $|\nabla \Phi (\xi)|\gtrsim |\xi|^{\alpha-1}$ where $\alpha\geq 1$ and $\gamma$ is a multi-index, which in particular involves the fractional Schr\"odinger operator $e^{-it(-\Delta)^{\frac{\alpha}{2}}}$; see \cite{lee2006pointwise} and \cite{cho2018note}.

Meanwhile further generalizations were established using geometric measure theory. Though Carleson's problem has an affirmative answer for a.e. convergence when $s \in [\frac{1}{4},\frac{1}{2}]$ for $n=1$, the divergence set (points $x\in \mathbb R$ where divergence occurs), which is of Lebesgue measure zero for such $s$, can still be big. \cite{barcelo2011dimension} shows that the divergence set is of Hausdorff dimension at most $1-2s$ for $s \in [\frac{1}{4},\frac{1}{2}]$. On the other hand, \cite{luca2017note} generalizes the necessity result of \cite{bourgain2016note} from the Lebesgue measure to the set of $\alpha$-dimensional non-negative measures $\mu$ on $\mathbb{R}^n$ for $n \geq 2$; here a non-negative Borel measure $\mu$ is $\alpha$-dimensional if $c_\alpha(\mu) = \sup \limits_{x\in \mathbb{R}^n,r>0} \frac{\mu(B(x,r))}{r^\alpha}<\infty$. It is shown that if $\alpha \in [\frac{3n+1}{4},n]$ and $\mu \left\{ x\in \mathbb{R}^n: e^{it\Delta}u_0 \xrightarrow[t\to 0]{}u_0 \:\text{fails}.\right\}=0$ for all $u_0 \in H^s(\mathbb{R})$, then $s \geq \frac{(n-1)(n-\alpha)}{2(n+1)}+\frac{n}{2(n+1)}$; since the Lebesgue measure on $\mathbb{R}^n$ is $n$-dimensional, the result of \cite{bourgain2016note} is recovered by letting $\alpha=n$. For recent results regarding the size of divergence set in higher dimensions, see \cite{du2018sharp}.

It offers some insight to view this convergence problem in the context of summation methods. These originated in the study of alternative ways of summing Fourier series  such as Abel or Riesz summability.   Summation methods for Fourier series or transforms, in modern terms, involve a family of operators $\phi(-t\Delta)$ (with $\phi$ a Borel-measurable or continuous function satisfying $\phi(0)=1$) forming an approximate identity as $t\rightarrow 0$.
Questions of convergence in this context translate into strong convergence (as $t\rightarrow 0$) of such operator families.  Abel summability corresponds to $\phi(x)=e^{-x}$, while other methods correspond to different choices of $\phi$ with $\phi(0)=1$. Our current (Schr\"odinger) problem chooses $\phi(x)=e^{ix}$, while the original result of Carleson for a.e.\ convergence of Fourier series \cite{carleson1966convergence} made the analogous statement for $\phi(x)=\frac{\sin(x)}{x}$.

The main purpose of this paper is to answer a variant of Carleson's problem, not for the free Schr\"odinger equation, but for the Schr\"odinger equation with a nonzero potential or nonlinearity. One motivation of this note comes from \cite{cowling1983pointwise} that whenever $|H|^\alpha u_0 \in L^2(X)$ where $X$ is a measure space and $H$ is some self-adjoint operator on $L^2(X)$ with $|H|$ given by the polar decomposition, we obtain $e^{-itH}u_0\xrightarrow[t\to 0]{}u_0$ a.e. if $\alpha\in (\frac{1}{2},\infty)$. Another motivation comes from \cite{sjogren2009boundary}, where given the following Cauchy problem, a.k.a. the quantum harmonic oscillator,
\begin{equation*}
    \begin{cases}
    i\partial_t u = -\partial_{xx}u + x^2u, (x,t) \in \mathbb{R}\times \mathbb{R}\\
    u(0)=u_0 \in H^s(\mathbb{R}),
    \end{cases}  
\end{equation*}
pointwise convergence to initial data holds for every $s \geq \frac{1}{4}$ and fails for $s<\frac{1}{4}$. Typically, a standard strategy in proving such positive result is to show that the Schr\"odinger maximal operator satisfies either a strong-type or weak-type estimate, from which pointwise convergence follows by a now-standard approximation argument. For the quantum harmonic oscillator, \cite{sjogren2009boundary} takes advantage of the closed, analytic expression for the fundamental solution associated with the quadratic Schr\"odinger propagator, also known as the Mehler kernel:
\begin{equation*}
       K_{it}(x,y) = (2\pi i \sin{2t})^{-n/2} e^{\frac{i}{2}(\cot{2t}\cdot |y-\frac{x}{\cos{2t}}|^2 - \tan{2t}\cdot|x|^2)}, \forall x,y \in \mathbb{R}^n, t \in \mathbb{R}\setminus \frac{\pi}{2}\mathbb{Z}. 
\end{equation*}

For a general potential, we have to work with analytic properties of the unitary group generated by the Hamiltonian $-\Delta+V$; note that the semigroup generated by this operator has been studied extensively, for example, by \cite{simon1982schrodinger}. In fact, an orbit of a square-integrable function generated by $e^{t\Delta}$, viewed as a spacetime function, solves the heat equation, and by exploiting the exponential decay of the corresponding Green's function, one can easily show pointwise convergence to initial data (the Green's function corresponding to $e^{it\Delta}$ has no such spatial decay). More generally, $\left\{e^{t\Delta}\right\}_{t>0}$ defines a holomorphic $C_0$-semigroup, and the strong convergence $e^{t\Delta}\xrightarrow[t\rightarrow 0]{}\mathds{I}$ is an example of standard Abel summability traditionally studied for Fourier series on an interval.  For complex $t\rightarrow 0$ such convergence occurs in a sector symmetric about the positive $t$ axis. However under the Wick rotation $t \mapsto it$, our sector of convergence is now symmetric about the imaginary $t$ axis, and our case of real $t\rightarrow 0$ constitutes a boundary case of the known region of Abel summability.  Therefore Abel summation is an insufficient tool to answer our problem. To this end, we summarise the main results of this paper:

\begin{theorem}\label{Mainresult}
Suppose $s\geq \frac{1}{4}$ and $V\in L^2(\mathbb{R})\cup \Big(W^{1,\infty}(\mathbb{R})\cap \bigcup\limits_{\rho\in [1,\infty)} L^\rho(\mathbb{R})\Big)$. Then the solutions to the linear Schr\"odinger equation converge a.e. to initial data in $H^s(\mathbb{R})$. On the other hand, if $s<\frac{1}{4}$ and $V\in L^2(\mathbb{R})$, then there exists a compactly supported initial data $f\in H^s(\mathbb{R})$ and a measurable set (of positive measure) $E_f \subseteq \mathbb{R}\setminus supp(f)$ such that $\varlimsup\limits_{t\rightarrow 0} |e^{it\partial_{xx}}f|\geq c>0$ on $E_f$.
\end{theorem}

It is clear from \cref{Mainresult} that if $V \in L^2(\mathbb{R})$, then the pointwise convergence of interest holds for $s \geq \frac{1}{4}$ sharply. Moreover the main theorem of \cite{carleson1980some} is contained in the previous statement by taking $V=0$. This class of potentials $V$ contains some well-studied examples in physics such as the finite square well.

The above results related to the linear Schr\"odinger equation are naturally related to corresponding non-linearizations (\cite{kenig1996quadratic},\cite{bejenaru2006sharp}), for which, perhaps as expected, the corresponding results hold.

\begin{theorem}\label{main2}
The solutions to quadratic nonlinear Schr\"odinger equation (qNLS) with nonlinearities 
\begin{center}
    $N_1(u,\overline{u})=u^2;\:N_2(u,\overline{u})=u\overline{u};\:N_3(u,\overline{u})=\overline{u}^2$  
\end{center} 
converge a.e. to initial data for $s\geq \frac{1}{4}$, $s>\frac{1}{4}$ and $s\geq \frac{1}{4}$, respectively. On the other hand, the convergence fails for qNLS with nonlinearities $N_1$ and $N_3$ in $H^s(\mathbb{R})$ with $s\in [0,\frac{1}{4})$.
\end{theorem}

We outline the organization of this article. In \cref{prelim}, useful notations are introduced. In \cref{positive1,positive2}, we prove a positive pointwise convergence result for the linear Schr\"odinger equation with potential using restricted Fourier space methods and Trotter-Kato product formula. In fact the class of potentials investigated does not include the quadratic case $V(x)=x^2$; our choice of potentials should be thought of as small perturbations to the free case $V=0$. In \cref{qnonlinearity}, the quadratic nonlinearities are treated. In \cref{negative2}, we prove the negative result that for $i\partial_t u = -\partial_{xx}u +Vu$, with an appropriate potential function, to exhibit pointwise convergence to initial data, it is necessary that $u_0 \in H^s(\mathbb{R})$ where $s\geq \frac{1}{4}$.

\section{Notation and Preliminaries.}\label{prelim}
The spaces $\mathscr{S}(\mathbb{R}^n)$ and $C^\infty_c(\mathbb{R}^n)$ denote the Schwartz class of rapidly decaying smooth functions and the set of smooth functions with compact support, respectively. We fix $\eta\in C^\infty_c(\mathbb{R})$ to be a smooth cutoff function that is identically one on $[-1,1]$ with a compact support in $[-2,2]$. The inhomogeneous and homogeneous differential operators are:
\begin{equation*}
\langle \xi \rangle = (1+|\xi|^2)^{1/2};\: \langle \nabla \rangle^s f = \mathcal{F}^{-1}(\langle \xi \rangle^s \Hat{f});\:| \nabla |^s f = \mathcal{F}^{-1}(| \xi |^s \Hat{f}),
\end{equation*}
and $\langle \partial_x\rangle, |\partial_x|,\langle \partial_t\rangle$ and $|\partial_t|$ are defined similarly.

The $L^2$-based Sobolev space and (dispersive) Sobolev spacess (also known as Fourier restriction space or Bourgain space in the literature) are:
\begin{equation*}
\begin{split}
H^s(\mathbb{R}^n) &= \langle \nabla \rangle^{-s} L^2(\mathbb{R}^n);\: \lVert f \rVert_{H^s} = \lVert \langle \nabla \rangle^s f \rVert_{L^2}\\
    X^{s,b} &= \left\{u:\lVert u \rVert_{X^{s,b}}<\infty\right\};\: \lVert u \rVert_{X^{s,b}} = \lVert \hat{u}(\xi,\tau) \langle \xi \rangle^{s}\langle \tau + \xi^2 \rangle^b \rVert_{L^2_{\xi,\tau}}. 
\end{split}
\end{equation*}

To do a local-in-time argument, where $t \in [-\delta,\delta]$ for some $\delta \in (0,1]$, we will need a restricted version of $X^{s,b}$ as well. We denote such a space by $X^{s,b}_\delta$ and its restricted norm as:
\begin{equation*}
       \lVert u \rVert_{X^{s,b}_\delta} = \inf \limits_{u=\Tilde{u}, t \in [-\delta,\delta]} \lVert \Tilde{u}\rVert_{X^{s,b}}. 
\end{equation*}

For $k\in\mathbb{N}$, the inhomogeneous and homogeneous $L^\infty$-based Sobolev spaces are:
\begin{equation*}
\begin{split}
W^{k,\infty} &= \left\{f: \lVert f \rVert_{W^{k,\infty}}<\infty\right\};\: \lVert f \rVert_{W^{k,\infty}} = \sum\limits_{|\alpha|\leq k} \lVert \partial^\alpha f \rVert_{L^\infty}\\
    \dot{W}^{k,\infty} &= \left\{ f: \lVert f \rVert_{\dot{W}^{k,\infty}}<\infty \right\};\: \lVert f \rVert_{\dot{W}^{k,\infty}} = \sum\limits_{|\alpha|= k} \lVert \partial^\alpha f \rVert_{L^\infty}.    
\end{split}
\end{equation*}
For $\delta>0$, the Banach space of continuous spacetime functions $u:[-\delta,\delta]\subseteq \mathbb{R} \rightarrow H^s(\mathbb{R}^n)$ are denoted by $C^0_tH^s_x([-\delta,\delta]\times\mathbb{R}^n)$ where $\lVert u \rVert_{C^0_tH^s_x} = \sup\limits_{t\in [-\delta,\delta]} \lVert u(t)\rVert_{H^s(\mathbb{R}^n)}$. Let $\mathcal{H} = p.v.(\frac{1}{x})$ be the Hilbert transform on $\mathbb{R}$. By Fourier analysis, $\mathcal{H}= \mathcal{F}^{-1}\Big(-i sgn(\xi) \Big)\mathcal{F}$ and hence defines a unitary operator on $L^2(\mathbb{R})$, and moreover $|\partial_x| = \partial_x \mathcal{H}$, again by considering their Fourier multipliers.

We say $A \lesssim B$ if $A$ is bounded above by $B$ multiplied by a universal constant, i.e., if there exists $C>0$ such that $A \leq CB$. Similarly, say $A \sim B$ if $A \lesssim B$ and $B \lesssim A$. For a measurable set $E\subseteq \mathbb{R}$, we let $|E|$ be the Lebesgue measure of $E$. We define $s+ = s+\epsilon$ for some universal $\epsilon <<1$; $s-$ is defined similarly. We assume $s\geq 0$ unless stated otherwise.

Lastly, some well-known properties of $X^{s,b}$ space and basic calculus facts are stated.

\begin{lemma}\label{l:32} Let $\delta \in (0,1]$. Then
\begin{enumerate}
    \item \cite{tao2006nonlinear}: For every $b > \frac{1}{2}$ and $\delta>0$, the following continuous embedding holds: $X^{s,b}_\delta \hookrightarrow C^0_tH^s_x([-\delta,\delta]\times\mathbb{R})$.
    \item \cite{tao2006nonlinear}: Linear estimate: $\lVert e^{it\partial_{xx}}f\rVert_{X^{s,b}_{\delta}}\lesssim_{s,b}\lVert f \rVert_{H^s}$ whenever the right-hand side is finite.
    \item \cite{tao2006nonlinear}: Let $b \in (\frac{1}{2},1]$ and $s \in \mathbb{R}$. Then, $\lVert \eta(t) \int_0^t e^{i(t-\tau)\partial_{xx}} F(\tau)d\tau \rVert_{X^{s,b}_\delta} \lesssim \lVert F \rVert_{X^{s,b-1}_\delta}$.
    \item \cite{taylor2007tools}: For $s \geq 0$ and $\frac{1}{2}=\frac{1}{p_1}+\frac{1}{q_1}=\frac{1}{p_2}+\frac{1}{q_2}$ where $p_1,q_2 \in (2,\infty]$, we have the following Leibniz rule for the $L^2$-based Sobolev space: $\lVert fg \rVert_{H^s} \lesssim \lVert f \rVert_{L^{p_1}} \lVert \langle \nabla \rangle^s g \rVert_{L^{q_1}}+\lVert \langle \nabla \rangle^s f \rVert_{L^{p_2}}\lVert g \rVert_{L^{q_2}}$
    \item \cite{erdougan2013smoothing}: If $\beta \geq \gamma \geq 0$ and $\beta + \gamma>1$, then
    \begin{equation*}
        \int \dfrac{dx}{\langle x-a_1 \rangle^\beta \langle x-a_2 \rangle^\gamma} \lesssim \langle a_1-a_2 \rangle^{-\gamma} \phi_\beta (a_1-a_2),
    \end{equation*}
    where\\
    \begin{center}
        $\phi_\beta (a) \sim \begin{cases}
    1,& \beta > 1\\
    \log (1+\langle a \rangle),              & \beta=1\\
    \langle a \rangle^{1-\beta}, & \beta<1
\end{cases}$
    \end{center}
    \item \cite{tao2006nonlinear}: Let $-\frac{1}{2}<b^\prime\leq b<\frac{1}{2}, s\in \mathbb{R}$ and $\delta\in (0,1]$. Then, $\lVert u\rVert_{X^{s,b^\prime}_{\delta}}\lesssim_{s,b^\prime,b}\delta^{b-b^\prime}\lVert u \rVert_{X^{s,b}_{\delta}}$.
    \item For all $a,b\in\mathbb{R}$, we have $\langle \tau-a \rangle \langle \tau - b \rangle \gtrsim \langle a-b \rangle$.
\end{enumerate}
\end{lemma}

\section{Linear Operator Estimates: Positive Results.}\label{positive1}
Let $H = -\partial_{xx}+V$ denote the Hamiltonian operator on  $\mathbb R$, where $V=V(x)$ is a real-valued time-independent multiplication operator. Note that $H$ is self-adjoint on $D(H)=D(-\partial_{xx})\cap D(V)$, if $V \in L^2(\mathbb{R})\cup  L^\infty(\mathbb{R})$, where $D(-\partial_{xx})= H^2(\mathbb{R})$ and $D(V)=\left\{f\in L^2(\mathbb{R}): Vf\in L^2(\mathbb{R})\right\}$; see \cite[Theorem 9.38]{hall2013quantum}. Therefore, $e^{-itH}$ gives a family of unitary actions on $L^2(\mathbb{R})$. It is of interest to ask whether known positive results for pointwise convergence of the free Schr\"odinger equation as $t\rightarrow 0$ can be recovered with an addition of a potential.

\begin{theorem}\label{mainresult}
Let $s\geq \frac{1}{4}$ and $1\leq \rho <\infty$, and suppose a time-independent potential $V$ satisfies the following hypothesis:
\begin{equation*}
    V \in L^2(\mathbb{R})\cup \Big(W^{1,\infty}(\mathbb{R})\cap L^\rho(\mathbb{R})\Big).
\end{equation*}
Then for all $u_0 \in H^s(\mathbb{R})$, $e^{-itH}u_0 \rightarrow u_0$ as $t\rightarrow 0$ almost everywhere with respect to Lebesgue measure. More precisely,
\begin{equation*}
    \left|\left\{x\in\mathbb{R}:\varlimsup\limits_{t\rightarrow 0} |e^{-itH}u_0-u_0|>0\right\}\right|=0.
\end{equation*}
\end{theorem}
\begin{remark}
By virtue of $V$ being time-independent, the conclusion holds in the limit when $t\rightarrow t_0$ for any $t_0 \in I\subseteq\mathbb{R}$, i.e., $e^{-itH}u_0\xrightarrow[t\rightarrow t_0]{}e^{-it_0 H}u_0$ a.e.
\end{remark} 
By Stone's theorem on a Hilbert space, a time-evolution operator for non-relativistic quantum mechanics is in one-to-one correspondence with a self-adjoint operator. However, the self-adjointness of $H$ generally fails on $H^s(\mathbb{R})$ for $s>0$, and therefore, $e^{-itH}$ defines a family of unitary operators on $H^s(\mathbb{R})$ only if $s=0$. In fact, it is not clear whether we have persistence of regularity for $e^{-itH}$ on $H^s(\mathbb{R})$ for $s>0$, and so this shall be proved. Some of these results are likely to be known; however the lemmas below contain some estimates that will be of use later. We begin with definitions (see \cite{tao2006nonlinear}).

\begin{definition}
For $\delta>0$, $u\in C^0_tH^s_x([-\delta,\delta],\mathbb{R})$ is a \textit{strong solution} of
\begin{equation}\label{maineq}
    \begin{cases}
    i\partial_t u &= -\partial_{xx}u + Vu, (x,t) \in \mathbb{R}\times [-\delta,\delta]\\
    u(0)&=u_0 \in H^s(\mathbb{R}),
    \end{cases}
\end{equation}
if $u$ satisfies the following Duhamel integral formula for all $t\in [-\delta,\delta]$:
\begin{equation*}
    u(t)=e^{it\partial_{xx}}u_0-i\int_0^t e^{i(t-t^\prime)\partial_{xx}}(Vu)(t^\prime)dt^\prime.
\end{equation*}
\end{definition} 

\begin{definition}
The Cauchy problem \cref{maineq} is \textit{well-posed} in $H^s(\mathbb{R})$ if for every $g\in H^s(\mathbb{R})$, there exists $\delta>0$, an open ball $B\subseteq H^s(\mathbb{R})$ containing $g$, and a subset $X\subseteq C^0_tH^s_x([-\delta,\delta]\times\mathbb{R})$ such that for every $u_0\in B$, there exists a unique strong solution $u\in X$ whose map $u_0\mapsto u$ is continuous. If $\delta>0$ can be arbitrarily large, then we say the well-posedness is \textit{global}.
\end{definition} 

\begin{remark}
For $u_0\in L^2(\mathbb{R})$, we claim that the notion of strong solution as in above, where we treat the potential term as a perturbation, coincides with that of an orbit generated by the unitary group. Though this seems intuitive, some care is needed if $V$ is not sufficiently regular. At least when $u_0\in D(H)$, $u(t)=e^{-itH}u_0$ satisfies the Duhamel integral formula for each $t$, which is an immediate consequence of the following product rule:
\begin{equation*}
   \partial_t \Big(e^{-it\partial_{xx}}e^{-itH}u_0 \Big) = -ie^{-it\partial_{xx}}\Big(Ve^{-itH}u_0\Big). 
\end{equation*}

For $u_0 \in L^2(\mathbb{R})\setminus D(H)$, let $u_0^{(n)}\rightarrow u_0$ as $n\rightarrow\infty$ where $u_0^{(n)}\in D(H)$. Then we have,
\begin{equation*}
    e^{-itH}u_0^{(n)}=e^{it\partial_{xx}}u_0^{(n)}-i\int_0^t e^{i(t-t^\prime)\partial_{xx}}(Ve^{-it^\prime H}u_0^{(n)})dt^\prime.
\end{equation*}

As $n\rightarrow\infty$, we have $e^{-itH}u_0^{(n)}\rightarrow e^{-itH}u_0$ and $e^{it\partial_{xx}}u_0^{(n)}\rightarrow e^{it\partial_{xx}}u_0$ by unitarity. We claim
\begin{equation*}
\int_0^t e^{i(t-t^\prime)\partial_{xx}}(Ve^{-it^\prime H}u_0^{(n)})dt^\prime\rightarrow \int_0^t e^{i(t-t^\prime)\partial_{xx}}(Ve^{-it^\prime H}u_0)dt^\prime,   
\end{equation*}
in $L^2(\mathbb{R})$ as $n\rightarrow\infty$. Firstly for $V\in L^\infty(\mathbb{R})$, we have
\begin{equation*}
\begin{split}
\left|\left|\int_0^t e^{i(t-t^\prime)\partial_{xx}}\Big(V(e^{-it^\prime H}u_0^{(n)}-e^{-it^\prime H}u_0)\Big)dt^\prime\right|\right|_{L^2}&\leq \int_0^t \lVert \Big(V(e^{-it^\prime H}u_0^{(n)}-e^{-it^\prime H}u_0)\Big)\rVert_{L^2}\\
    &\leq \int_0^t \lVert V \rVert_{L^\infty}\lVert u_0^{(n)}-u_0\rVert_{L^2}dt^\prime \rightarrow 0,    
\end{split}
\end{equation*}
where the last inequality is by the H\"older's inequality.

Secondly for $V\in L^2(\mathbb{R})$, we apply the following form of inhomogeneous Strichartz estimate (see \cite[Theorem 2.3]{tao2006nonlinear}): 
\begin{equation*}
    \left|\left| \int_0^t e^{i(t-t^\prime)\partial_{xx}}F(t^\prime)dt^\prime\right|\right|_{L^\infty_t L^2_x}\lesssim \lVert F \rVert_{L^{\frac{4}{3}}_t L^1_x}.
\end{equation*}

Let $\chi$ be a characteristic function on $t^\prime \in [0,T]$ where $0<t<T$. Then we have,
\begin{equation*}
\begin{split}
\left|\left| \int_0^t e^{i(t-t^\prime)}\Big(V e^{-it^\prime H}(u_0^{(n)}-u_0)\Big)dt^\prime\right|\right|_{L^2_x}&=\left|\left| \int_0^t e^{i(t-t^\prime)}\Big(V\chi(t^\prime) e^{-it^\prime H}(u_0^{(n)}-u_0)\Big)dt^\prime\right|\right|_{L^2_x}\\
    &\leq \left|\left| \int_0^t e^{i(t-t^\prime)}\Big(V\chi(t^\prime) e^{-it^\prime H}(u_0^{(n)}-u_0)\Big)dt^\prime\right|\right|_{L^\infty_t L^2_x}\\
    &\lesssim \lVert V\chi e^{-itH}(u_0^{(n)}-u_0) \rVert_{L_t^{\frac{4}{3}}L^1_x}\lesssim \lVert V \rVert_{L^2} \lVert u_0^{(n)}-u_0\rVert_{L^2}\rightarrow 0.    
\end{split}
\end{equation*}
\end{remark} 

The key idea of our proof is the Bourgain space estimate of the potential term.

\begin{lemma}\label{l:33}
Let $s\in [0,\frac{3}{4})$. Then, there exists $b \in (\frac{1}{2},1]$, $\gamma \in [0,\frac{1}{2})$ and $a\in [0,\frac{1}{2})$ that satisfy
\begin{equation*}
    \frac{s+a}{2}\leq \gamma <\min(\frac{s}{2}+\frac{1}{4},\frac{1}{2});\:\max(\frac{s}{2}+\frac{1}{4},\frac{1}{2})< b <1-\gamma.
\end{equation*}
Furthermore for every such $(s,b,\gamma,a)$, we have $\lVert Vu \rVert_{X^{s+a,-\gamma}_\delta}\lesssim_{s,b,\gamma,a}\lVert V \rVert_{L^2}\lVert u \rVert_{X^{s,b}_\delta}$.
\end{lemma}
\begin{proof}[proof of \cref{l:33}]
The first statement is a  straightforward algebra exercise. As for the second, it suffices to prove the statement neglecting the $\delta$-dependence, for if $\Tilde{u} = u$ on $t\in[-\delta,\delta]$, we have
\begin{equation*}
    \lVert Vu \rVert_{X^{s,-\gamma}_\delta}\leq \lVert \eta(\frac{t}{\delta})V\Tilde{u}\rVert_{X^{s,-\gamma}}\lesssim_\eta \lVert V\Tilde{u}\rVert_{X^{s,-\gamma}}\lesssim \lVert V \rVert_{L^2}\lVert \Tilde{u}\rVert_{X^{s,b}}.
\end{equation*}
Taking infimum over $\Tilde{u}$, we derive the desired result. We argue as in the proof of \cite[Proposition 1]{erdogan2013talbot}.
Define
\begin{equation*}
    F(\xi) = |\hat{V}(\xi)|;\: G(\xi,\tau) = \langle \xi\rangle^{s+a} \langle \tau + \xi^2 \rangle^b |\hat{u}(\xi,\tau)|,
\end{equation*}
and
\begin{equation*}
    W(\xi, \tau, \xi_1) = \dfrac{\langle \xi \rangle^{s+a} \langle \tau + \xi^2 \rangle^{-\gamma}}{\langle \xi_1 \rangle^{s} \langle \tau + \xi_1^2 \rangle^b}.
\end{equation*}
Noting that $\mathcal{F}[Vu](\xi,\tau) = \int \hat{V}(\xi-\xi_1)\hat{u}(\xi_1,\tau)d\xi_1$, we have $ \lVert Vu \rVert_{X^{s,-\gamma}}^2$
\begin{equation*}
\begin{split}
    &= \left|\left| \int \langle \xi \rangle^{s+a} \langle \tau + \xi^2 \rangle^{-\gamma}\hat{V}(\xi-\xi_1)\hat{u}(\xi_1,\tau)d\xi_1\right|\right|_{L^2_{\xi,\tau}}^2
    \leq \left|\left| \int \dfrac{\langle \xi \rangle^{s+a} \langle \tau + \xi^2 \rangle^{-\gamma}}{\langle \xi_1 \rangle^{s} \langle \tau + \xi_1^2 \rangle^b} F(\xi-\xi_1)G(\xi_1,\tau)d\xi_1\right|\right|_{L^2_{\xi,\tau}}^2\\
    &\leq \left|\left| \Big(\int W^2 d \xi_1 \Big)^{1/2} \Big(\int F(\xi-\xi_1)^2G(\xi_1,\tau)^2d\xi_1 \Big)^{1/2}\right|\right|_{L^2_{\xi,\tau}}^2
    = \left|\left| \int W^2 d\xi_1 \cdot\int F(\xi-\xi_1)^2 G(\xi_1,\tau)^2d\xi_1\right|\right|_{L^1_{\xi,\tau}}\\
    &\leq \left|\left| \int W^2 d\xi_1\right|\right|_{L^\infty_{\xi,\tau}} \cdot\lVert F^2 \ast_{\xi_1} G^2\rVert_{L^1_{\xi,\tau}}\leq \left|\left| \int W^2 d\xi_1\right|\right|_{L^\infty_{\xi,\tau}} \cdot \lVert V \rVert_{L^2}^2 \lVert u\rVert_{X^{s,b}}^2,
\end{split}    
\end{equation*}
where the second inequality is due to the Cauchy-Schwarz inequality, the third by the H\"older's inequality and the fourth by the Young's inequality. It remains to prove that $\lVert \int W^2 d\xi_1\rVert_{L^\infty_{\xi,\tau}}$ is finite. Changing variable $z=\xi_1^2$,
\begin{equation*}
    \int W^2 d\xi_1 \simeq \langle \xi \rangle^{2s+2a} \langle \tau + \xi^2 \rangle^{-2\gamma} \int \frac{d\xi_1}{\langle \xi_1^2 \rangle^{s} \langle \tau+ \xi_1^2 \rangle^{2b}}\simeq \langle \xi \rangle^{2s+2a} \langle \tau + \xi^2 \rangle^{-2\gamma} \int_0^\infty \frac{dz}{\langle z \rangle^s\langle z+\tau \rangle^{2b}z^{1/2}}.
\end{equation*}
Note that $\sup\limits_{\tau\in\mathbb{R}}$ can be replaced by $\sup\limits_{|\tau|>1}$ without loss of generality, for if $|\tau|\leq 1$, then $\langle \tau+\xi^2\rangle^{-2\gamma}\leq \langle \xi^2-1\rangle^{-2\gamma}$ for $|\xi|\geq 1$, and 
\begin{equation*}
    \sup\limits_{|\tau|\leq 1}\int_0^\infty \frac{dz}{\langle z \rangle^s\langle z+\tau \rangle^{2b}z^{1/2}}<\infty.
\end{equation*}

Hence $\sup\limits_{|\xi|\geq 1,|\tau|\leq 1}\int W^2 d\xi_1\lesssim \sup\limits_{|\xi|\geq 1}\langle \xi \rangle^{2s+2a-4\gamma}<\infty$, whereas $\sup\limits_{|\xi|\leq 1,|\tau|\leq 1}\int W^2d\xi_1<\infty$ follows from extreme value theorem. Now suppose $|\tau|>1$. Then,
\begin{equation*}
    \begin{split}
        \int_0^\infty \frac{dz}{\langle z \rangle^s\langle z+\tau \rangle^{2b}z^{1/2}}&=\int_0^\frac{1}{2} \frac{dz}{\langle z \rangle^s\langle z+\tau \rangle^{2b}z^{1/2}}+\int_{\frac{1}{2}}^\infty \frac{dz}{\langle z \rangle^s\langle z+\tau \rangle^{2b}z^{1/2}}\\
     &\lesssim \langle \tau \rangle^{-2b} + \langle \tau \rangle^{-(s+\frac{1}{2})}\lesssim \langle \tau \rangle^{-(s+\frac{1}{2})},
    \end{split}
\end{equation*}
since $s+\frac{1}{2}<2b$. Moreover since $2\gamma < s+\frac{1}{2}$,
\begin{equation*}
\sup\limits_{\xi\in\mathbb{R},|\tau|>1}\langle \xi \rangle^{2s+2a} \langle \tau + \xi^2 \rangle^{-2\gamma}\langle \tau \rangle^{-(s+\frac{1}{2})}\lesssim \sup\limits_{\xi\in\mathbb{R}}\langle \xi \rangle^{2s+2a-4\gamma}<\infty.
\end{equation*}
\end{proof}

When $V \in L^2(\mathbb{R})$, the following Bourgain space estimate is obtained with ease via Fourier analysis.

\begin{lemma}\label{l:31}
Suppose $V \in L^2(\mathbb{R})$. The Cauchy problem \cref{maineq} is globally well-posed in $H^s(\mathbb{R})$ for $s \in [0,\frac{3}{4})$. In particular if $u$ is the strong solution with the initial data $u_0 \in H^s(\mathbb{R})$, then there exists $\delta = \delta(\lVert V \rVert_{L^2})>0$ such that $\lVert u \rVert_{X^{s,b}_\delta}\lesssim_{s,b,\lVert V \rVert_{L^2}} \lVert u_0\rVert_{H^s}$ for some $b \in (\frac{1}{2},1]$.
\end{lemma} 

\begin{proof}[proof of \cref{l:31}]
Assuming that \cref{l:33} holds, let $s,b,\gamma$ be as in \cref{l:33}, $\delta\in (0,1]$ and fix $C>0$ that satisfies $\lVert e^{it\partial_{xx}}f\rVert_{X^{s,b}_\delta}\leq C\lVert f \rVert_{H^s}$ for all $f\in H^s(\mathbb{R})$ by \cref{l:32}. Let 
\begin{equation*}
X = \left\{u\in X^{s,b}_\delta:\lVert u \rVert_{X^{s,b}_\delta}\leq 2C\lVert u_0\rVert_{H^s}\right\}.    
\end{equation*}

Define $\Gamma u = e^{it\partial_{xx}}u_0 - i \int_0^t e^{i(t-t^\prime)\partial_{xx}}(Vu)(t^\prime)dt^\prime$. Then by \cref{l:32,l:33} we have,
\begin{equation*}
\begin{split}
    \lVert \Gamma u \rVert_{X^{s,b}_\delta}&\lesssim \lVert u_0 \rVert_{H^s} + \lVert Vu \rVert_{X^{s,b-1}_\delta}\lesssim \lVert u_0 \rVert_{H^s}+\delta^{1-(b+\gamma)}\lVert Vu \rVert_{X^{s,-\gamma}_\delta}\\
    &\lesssim \lVert u_0 \rVert_{H^s}+\delta^{1-(b+\gamma)}\lVert V\rVert_{L^2}\lVert u \rVert_{X^{s,b}_\delta}.\\
    \Rightarrow \lVert \Gamma u \rVert_{X^{s,b}_\delta}&\leq C \lVert u_0\rVert_{H^s}+\Tilde{C}C\delta^{1-(b+\gamma)}\lVert V \rVert_{L^2}\lVert u_0\rVert_{H^s}.
\end{split}
\end{equation*}

By choosing $\delta \leq (\Tilde{C}\lVert V \rVert_{L^2})^{-\frac{1}{1-(b+\gamma)}}$, it is shown that $\Gamma:X\rightarrow X$. Similarly, we obtain
\begin{equation*}
    \lVert \Gamma u-\Gamma v \rVert_{X^{s,b}_\delta}\leq C_0 \delta^{1-(b+\gamma)}\lVert V \rVert_{L^2}\lVert u-v \rVert_{X^{s,b}_\delta},
\end{equation*}
from which it is shown that $\Gamma$ is a contraction map by shrinking $\delta>0$ if necessary, and the resulting unique fixed point is the desired strong solution. Since the time step only depends on the norm of $V$, this local result can be iterated infinitely many times, and hence our solution is global in time.

Continuous dependence on initial data follows similarly, for if $T>0$, $u_0^{(n)}\rightarrow u_0$ in $H^s(\mathbb{R})$ and $u^{(n)},u$ denote the strong solution corresponding to $u_0^{(n)}, u_0$, respectively, then for $t\leq T$,
\begin{equation*}
\begin{split}
    \lVert u^{(n)}(t)-u(t) \rVert_{H^s(\mathbb{R})}&\leq \lVert u_0^{(n)}-u_0\rVert_{H^s}+\left|\left| \int_0^t e^{i(t-t^\prime)\partial_{xx}}\Big(V(u^{(n)}-u)\Big)(t^\prime)dt^\prime\right|\right|_{H^s}\\
    &\lesssim \lVert u_0^{(n)}-u_0\rVert_{H^s}+T^{1-(b+\gamma)}\lVert V \rVert_{L^2}\lVert u^{(n)}-u\rVert_{X^{s,b}_T}\\
    &\lesssim \lVert u_0^{(n)}-u_0\rVert_{H^s}+T^{1-(b+\gamma)}\lVert V \rVert_{L^2}\lVert u^{(n)}_0-u_0\rVert_{H^s},
\end{split}
\end{equation*}
where the implicit constant may depend on $T$. Taking $\sup\limits_{t\in [0,T]}$ both sides and taking $n\rightarrow \infty$, we obtain the desired result.
\end{proof}

So far, the dispersive estimate of $e^{it\partial_{xx}}$ was used to control the Duhamel contribution of $V$. Now we directly study the dispersive estimate of $e^{-itH}$. If $V$ and $\partial_{xx}$ commute, then
\begin{equation}\label{commute}
e^{-itH}=e^{it\partial_{xx}}e^{-itV},
\end{equation}
and therefore, the operator $e^{-itH}$ would obey the same maximal operator estimate of $e^{it\partial_{xx}}$ as in \cite{carleson1980some}, and our problem would be trivial. Generally the exponential map does not take addition into multiplication. If $t$ is small, however, it is feasible to believe that \cref{commute} holds approximately, and the following lemma quantifies this intuition: 

\begin{lemma}\cite[Theorem 8.30]{simon1980methods}\label{l:41}
Let $A$ and $B$ be self-adjoint operators on a Hilbert space $\mathscr{H}$. If $A+B$ is self-adjoint on $D(A) \cap D(B)$, then
\begin{equation*}
\lim\limits_{n\rightarrow \infty} (e^{i\frac{t}{n}A}e^{i\frac{t}{n}B})^n\phi = e^{it(A+B)}\phi, \forall \phi \in \mathscr{H}.
\end{equation*}
\end{lemma}

We apply this \textit{Trotter-Kato product formula} to obtain persistence of regularity when the derivative of $V$ is bounded.

\begin{lemma}\label{l:42}
Suppose $t\in \mathbb{R}$ and $s\in[0,1]$. If $\lVert V \rVert_{\dot{W}^{1,\infty}(\mathbb{R}^n)}<\infty$, then we obtain
\begin{equation*}
\lVert e^{-itH}f\rVert_{H^s}\leq e^{st\sqrt{n}\lVert V \rVert_{\dot{W}^{1,\infty}}}\lVert f \rVert_{H^s}, \forall f\in C^\infty_c(\mathbb{R}^n).    
\end{equation*}
\end{lemma}

\begin{proof}[proof of \cref{l:42}]
We first show $\lVert e^{-itV}\rVert_{H^1 \rightarrow H^1}\leq 1+ t\sqrt{n}\lVert V \rVert_{\dot{W}^{1,\infty}}, \forall t \in \mathbb{R}$. Let $f \in C_c^\infty(\mathbb{R}^n)$. Then we have\\
\begin{equation*}
\begin{split}
    \lVert e^{-itV} f \rVert_{H^1}^2 &= \lVert e^{-itV} f \rVert_{L^2}^2 + \sum\limits_{j=1}^n \lVert \partial_j (e^{-itV}f)\rVert_{L^2}^2\\
    &\leq \lVert f \rVert_{L^2}^2 + \sum\limits_{j=1}^n (t\lVert \partial_j V \cdot f\rVert_{L^2}+\lVert \partial_j f \rVert_{L^2})^2\leq (1+t\sqrt{n}\lVert V \rVert_{\dot{W}^{1,\infty}})^2\lVert f \rVert_{H^1}^2.
\end{split}
\end{equation*}
Hence, the best constant $C(t) \leq 1+t\sqrt{n}\lVert V \rVert_{\dot{W}^{1,\infty}},\forall t \in \mathbb{R}$.\\

Let $\phi = e^{-itH}u_0$ for a fixed $u_0 \in H^1(\mathbb{R}^n)$ and $\phi_m = (e^{-i\frac{t}{m}V }e^{i \frac{t}{m}\partial_{xx}})^m u_0$. Then, $\phi_m \rightarrow \phi$ in $L^2(\mathbb{R}^n)$ by \cref{l:41}. By the estimate on $\lVert e^{-itV} \rVert_{H^1 \rightarrow H^1}$, we obtain $\lVert e^{-itV} e^{it\partial_{xx}} f\rVert_{H^1} \leq (1+t\sqrt{n}\lVert V \rVert_{\dot{W}^{1,\infty}})\lVert f \rVert_{H^1}$. Then we have
\begin{equation*}
    \lVert \phi_m \rVert_{H^1} \leq (1+\frac{\sqrt{n}\lVert V \rVert_{\dot{W}^{1,\infty}}t}{m})^m\lVert u_0 \rVert_{H^1}.
\end{equation*}
Hence for $t \in [0,T]$ for $T>0$, we have a bounded sequence $\left\{ \phi_m \right\}_m \subset H^1(\mathbb{R}^n)$, a reflexive Banach space. Then, there exists a weakly convergent subsequence $\left\{\phi_{m_k}\right\}_k$ where $\phi_{m_k} \rightharpoonup \Tilde{\phi} \in H^1(\mathbb{R}^n)$. Since $H^1(\mathbb{R}^n)\hookrightarrow L^2(\mathbb{R}^n)$, $\phi_{m_k} \rightharpoonup \Tilde{\phi}$ in $L^2(\mathbb{R}^n)$ and since $\phi_m \rightarrow \phi$ in $L^2(\mathbb{R}^n)$, the convergence holds in the weak topology, and by the uniqueness of weak-limit in Banach space, $\phi = \Tilde{\phi}$; in particular, $e^{-itH}u_0 \in H^1(\mathbb{R}^n)$. Since norm is lower semicontinuous with respect to weak topology,
\begin{equation*}
    \lVert e^{-itH}u_0 \rVert_{H^1} \leq \varliminf \limits_{k \rightarrow \infty} \lVert \phi_{m_k}\rVert_{H^1} \leq e^{t\sqrt{n}\lVert V \rVert_{\dot{W}^{1,\infty}} }\lVert u_0\rVert_{H^1}.
\end{equation*}
Since the bound above holds for all $t \in [0,T]$ uniformly in $T$, it holds for all $t \in \mathbb{R}$. Then by complex interpolation, it follows that for $s \in [0,1]$,
\begin{equation*}
    \lVert e^{-itH}f \rVert_{H^s}\leq e^{st\sqrt{n}\lVert V \rVert_{\dot{W}^{1,\infty}} }\lVert f\rVert_{H^s}, \forall f \in H^s(\mathbb{R}^n).
\end{equation*}
\end{proof}

Recall the estimate $\lVert e^{-itH}u_0 \rVert_{X^{s,b}_\delta}\lesssim \lVert u_0\rVert_{H^s}$ for $s\in [\frac{1}{4},\frac{1}{2}]$ and $b=\frac{1}{2}+$ was obtained in \cref{l:31}. A similar estimate via a different approach - the Trotter-Kato product formula and the fractional Gagliardo-Nirenberg interpolation - is obtained in \cref{l:43}. Let $s\in (0,1)$ and  $\max(\frac{1-s}{s},2)<\rho<\infty$. Then for every $f\in C^\infty_c(\mathbb{R})$\footnote{For a more general statement, see \cite[Theorem 1]{brezis2018gagliardo}.},
\begin{equation}\label{gagliardo}
    \lVert f \rVert_{W^{s,\frac{\rho}{1-s}}}\lesssim \lVert f \rVert_{L^\rho}^{1-s}\lVert f \rVert_{W^{1,\infty}}^s.
\end{equation}
\begin{lemma}\label{l:43}
Let $s\in (0,1),b \in (\frac{1}{2},1]$ and $\delta \in (0,1]$. Then, $\lVert e^{-itH} f \rVert_{X^{s,b}_\delta} \lesssim_{s,b,\eta} \lVert f \rVert_{H^s}, \forall f \in H^s(\mathbb{R})$.
\end{lemma}

\begin{proof}[proof of \cref{l:43}]
Let $F(\xi,\tau) = \mathcal{F}[\eta(\cdot)e^{-i\cdot H}f](\xi,\tau)$ and $\Tilde{F}(\xi,t) = \mathcal{F}^{-1}_{\tau}F$ where $\mathcal{F}^{-1}_{\tau}$ is the inverse Fourier transform in $\tau$ variable, and let $\mathcal{F}_x$ be defined similarly. Then we obtain
\begin{equation*}
\begin{split}
    \lVert e^{-itH} f \rVert_{X^{s,b}_\delta}&\leq
    \lVert \eta(t) e^{-itH}f \rVert_{X^{s,b}}= \lVert \langle \xi \rangle^s \langle \tau + \xi^2 \rangle^b F(\xi,\tau) \rVert_{L^2_{\xi,\tau}}\\
    &=\lVert \langle \xi \rangle^s \lVert e^{it\xi^2} \Tilde{F}(\xi,t)\rVert_{H^b_t} \rVert_{L^2_\xi}\leq \lVert \langle \xi \rangle^s \lVert e^{it\xi^2} \Tilde{F}(\xi,t)\rVert_{H^1_t} \rVert_{L^2_\xi}\\
    &\lesssim_{s} \lVert \langle \xi \rangle^s \Tilde{F}(\xi,t) \rVert_{L^2_{\xi,t}} + \lVert  \langle \xi \rangle^s\Big(|\partial_t| (e^{it\xi^2} \Tilde{F}(\xi,t)) \Big) \rVert_{L^2_{\xi,t}}.
\end{split}
\end{equation*}
For the first term, integrate in $\xi$ variable first using Plancherel's theorem, followed by the estimate for the operator norm $\lVert e^{-itH}\rVert_{H^s \rightarrow H^s}$ and followed by the $t$-integral as follows:
\begin{equation*}
    \lVert \langle \xi \rangle^s \Tilde{F}(\xi,t) \rVert_{L^2_{\xi,t}} =\lVert \eta(t)\cdot\lVert e^{-itH} f  \rVert_{H^s_x}\rVert_{L^2_t}\leq \lVert \eta(t) e^{s\lVert V^\prime \rVert_{L^\infty}t} \rVert_{L^2_t} \cdot \lVert f \rVert_{H^s}\lesssim_\eta \lVert f \rVert_{H^s}.
\end{equation*}
As for the second term,
\begin{equation*}
\begin{split}
    \lVert  \langle \xi \rangle^s\Big(|\partial_t| (e^{it\xi^2} \Tilde{F}(\xi,t)) \Big) \rVert_{L^2_{\xi,t}} &= \lVert \langle \xi \rangle^s \cdot \lVert |\partial_t| (e^{it\xi^2}\Tilde{F}(\xi,t)) \rVert_{L^2_t} \rVert_{L^2_\xi}=\lVert \langle \xi \rangle^s \cdot \lVert \mathcal{H} \partial_t e^{it\xi^2} \Tilde{F}(\xi,t) \rVert_{L^2_t}\rVert_{L^2_\xi}\\
    &=\lVert \langle \xi \rangle^s \cdot \lVert \partial_t \Big(\eta(t)e^{it\xi^2}\mathcal{F}_x[e^{-itH}f](\xi,t) \Big)\rVert_{L^2_t} \rVert_{L^2_\xi}\\
    &\leq \lVert \langle \xi \rangle^s \cdot \lVert \partial_t \eta \cdot e^{it\xi^2} \mathcal{F}_x[e^{-itH}f] \rVert_{L^2_t} \rVert_{L^2_\xi} + \lVert \langle \xi \rangle^s \cdot \lVert \eta(t) \partial_t \Big(e^{it\xi^2}\mathcal{F}_x[e^{-itH}f] \Big) \rVert_{L^2_t} \rVert_{L^2_\xi}.
\end{split}    
\end{equation*}
For the first term, switching the order of integration and recalling that the family $e^{it\partial_{xx}}$ is unitary on $H^s(\mathbb{R})$,
\begin{equation*}
\begin{split}
    \lVert \langle \xi \rangle^s \cdot \lVert \partial_t \eta \cdot e^{it\xi^2} \mathcal{F}_x[e^{-itH}f] \rVert_{L^2_t} \rVert_{L^2_\xi} &= \lVert \partial_t \eta \cdot \lVert e^{-it\partial_{xx}}e^{-itH}f \rVert_{H^s_x} \rVert_{L^2_t}\\
    &= \lVert \partial_t \eta \cdot\lVert e^{-itH}f \rVert_{H^s_x} \rVert_{L^2_t}\leq \lVert \partial_t \eta \cdot e^{s\lVert V^\prime \rVert_{L^\infty}t}\rVert_{L^2_t}\cdot \lVert f \rVert_{H^s}\lesssim_\eta \lVert f \rVert_{H^s}.
\end{split}    
\end{equation*}
For the second term, use product rule in $t$ to obtain
\begin{equation*}
    \lVert \langle \xi \rangle^s \cdot \lVert \eta(t) \partial_t \Big(e^{it\xi^2}\mathcal{F}_x[e^{-itH}f] \Big) \rVert_{L^2_t} \rVert_{L^2_\xi} = \lVert \eta(t)\cdot \lVert \partial_t \Big(e^{-it\partial_{xx}}e^{-itH}f \Big) \rVert_{H^s_x} \rVert_{L^2_t}= \lVert \eta(t) \cdot \lVert V e^{-itH}f \rVert_{H^s_x} \rVert_{L^2_t}.
\end{equation*}
where the second equality follows from $\partial_t \Big(e^{-it\partial_{xx}}e^{-itH}f \Big) = -ie^{-it\partial_{xx}}\Big(Ve^{-itH}f\Big)$. Then with $q \in (2,\infty)$ defined as follows,
\begin{equation*}
    \frac{1}{q}+\frac{1-s}{\rho}=\frac{1}{2},
\end{equation*}
apply the following particular form of Leibniz rule for Sobolev space to obtain\footnote{Unfortunately, the Leibniz rule generally fails when the $L^\infty$ norm is applied to the Bessel potential term. Had this been true, the decay condition on $V$ could have been removed.}
\begin{equation*}
\begin{split}
    \lVert V e^{-itH}f \rVert_{H^s} &\lesssim_s \lVert V \rVert_{L^\infty} \lVert e^{-itH}f \rVert_{H^s} +  \lVert V \rVert_{W^{s,\frac{\rho}{1-s}}} \lVert e^{-itH}f \rVert_{L^{q}}\\
    &\lesssim \Big(\lVert V \rVert_{L^\infty}+ \lVert V \rVert_{W^{s,\frac{\rho}{1-s}}}\Big)\lVert e^{-itH}f \rVert_{H^s} \leq (\lVert V \rVert_{L^\infty}+\lVert V \rVert_{W^{s,\frac{\rho}{1-s}}})e^{s \lVert V^\prime \rVert_{L^\infty} t} \lVert f \rVert_{H^s}.
\end{split}
\end{equation*}

Since the first factor of the RHS is finite by \cref{gagliardo}, the proof is complete by integrating the upper bound in $t$ against the smooth bump $\eta$.
\end{proof}

\begin{remark}
In the case of $V=0$, \cref{l:43} reduces to \cref{l:32} where the proof heavily depends on the fact that the time-evolution operator defines a Fourier multiplier. However, if $V$ is not identically zero, then the linear action by $e^{-itH}$ defines a Fourier integral operator. The linear estimate as above, therefore, is not entirely obvious for $e^{-itH}$.
\end{remark}

\section{Proof of \cref{mainresult}.}\label{positive2}

\begin{proof}[proof of \cref{mainresult}]
For initial data in $H^s(\mathbb{R})$ for $s>\frac{1}{2}$, the solution for each $t\in\mathbb{R}$ can be identified with a continuous function by Sobolev embedding, and therefore, the conclusion follows immediately. Let $V \in L^2(\mathbb{R})$. Suppose $u_0\in H^s(\mathbb{R})$ for $s\in [\frac{1}{4},\frac{1}{2}]$. We have
\begin{equation*}
\begin{split}
\left|\left\{x:\varlimsup\limits_{t\rightarrow 0}|e^{-itH}u_0 - u_0|>0\right\}\right| &\leq \left|\left\{x:\varlimsup\limits_{t\rightarrow 0}|e^{it\partial_{xx}}u_0 - u_0|>0\right\}\right|+\left|\left\{x:\varlimsup\limits_{t\rightarrow 0}\left|\int_0^t e^{i(t-t^\prime)\partial_{xx}}\Big(Vu\Big)(t^\prime)dt^\prime\right|>0\right\}\right|\\
    &=\left|\left\{x:\varlimsup\limits_{t\rightarrow 0}\left|\int_0^t e^{i(t-t^\prime)\partial_{xx}}\Big(Vu\Big)(t^\prime)dt^\prime\right|>0\right\}\right|,
\end{split}
\end{equation*}
where the equality holds due to \cite{carleson1980some}. Then there exists $a>0,b>\frac{1}{2}$ and $\gamma<\frac{1}{2}$ such that
\begin{equation*}
 s+a>\frac{1}{2};\: s+a\leq 2\gamma;\: b-1<-\gamma,
\end{equation*}
and for $\delta \in (0,1]$:
\begin{equation*}
    \lVert Vu\rVert_{X^{s+a,-\gamma}_\delta}\lesssim \lVert V\rVert_{L^2}\lVert u \rVert_{X^{s,b}_\delta}.
\end{equation*}

By \cref{l:32,l:33},
\begin{equation*}
\begin{split}
\left|\left| \int_0^t e^{i(t-t^\prime)\partial_{xx}}\Big(Vu\Big)(t^\prime)dt^\prime \right|\right|_{C^0_tH^{s+a}_x([-\delta,\delta]\times\mathbb{R})}&\lesssim \left|\left| \int_0^t e^{i(t-t^\prime)\partial_{xx}}\Big(Vu\Big)(t^\prime)dt^\prime\right|\right|_{X^{s+a,b}_\delta}\lesssim \lVert Vu\rVert_{X^{s+a,b-1}_\delta}\\
    &\lesssim \lVert V\rVert_{L^2}\lVert u \rVert_{X^{s,b}_\delta}\lesssim \lVert V \rVert_{L^2}\lVert u_0\rVert_{H^s}<\infty.
\end{split}
\end{equation*}

Hence another application of Sobolev embedding implies
\begin{equation*}
    \left|\left\{x:\varlimsup\limits_{t\rightarrow 0}\left|\int_0^t e^{i(t-t^\prime)\partial_{xx}}\Big(Vu\Big)(t^\prime)dt^\prime\right|>0\right\}\right|=0.
\end{equation*}

Now let $V\in W^{1,\infty}\cap L^\rho$. Fix an open cover $\left\{(\frac{k}{2},\frac{k}{2}+1)\right\}_{k\in\mathbb{Z}}$ of $\mathbb{R}$ and let $\left\{\psi_k\right\}_k$ be a smooth partition of unity subordinate to the open cover. Then, $V = \sum\limits_{k\in\mathbb{Z}} V_k$ where $V_k = V \psi_k\in L^2(\mathbb{R})$. For $u_0\in H^s(\mathbb{R})$ where $s\in [\frac{1}{4},\frac{1}{2}]$, we have
\begin{equation*}
\begin{split}
&\left|\left\{x:\varlimsup\limits_{t\rightarrow 0}|e^{-itH}u_0 - u_0|>0\right\}\right|\\ &\leq \left|\left\{x:\varlimsup\limits_{t\rightarrow 0}|e^{it\partial_{xx}}u_0 - u_0|>0\right\}\right|
    +\sum\limits_{k\in\mathbb{Z}}\left|\left\{x:\varlimsup\limits_{t\rightarrow 0}\left|\int_0^t e^{i(t-t^\prime)\partial_{xx}}\Big(V_ke^{-it^\prime H}u_0\Big)dt^\prime\right|>0\right\}\right|\\
    &=\sum\limits_{k\in\mathbb{Z}}\left|\left\{x:\varlimsup\limits_{t\rightarrow 0}\left|\int_0^t e^{i(t-t^\prime)\partial_{xx}}\Big(V_ke^{-it^\prime H}u_0\Big)dt^\prime\right|>0\right\}\right|.    
\end{split}
\end{equation*}
As before for $\delta \in (0,1]$ and $b=\frac{1}{2}+$, we obtain
\begin{equation*}
    \left|\left| \int_0^t e^{i(t-t^\prime)\partial_{xx}}\Big(V_ke^{-it^\prime H}u_0\Big)(t^\prime)dt^\prime \right|\right|_{C^0_tH^{\frac{1}{2}+}_x([-\delta,\delta]\times\mathbb{R})}\lesssim \lVert V_k\rVert_{L^2}\lVert e^{-it H}u_0 \rVert_{X^{\frac{1}{4},b}_\delta}\lesssim \lVert V_k \rVert_{L^2}\lVert u_0\rVert_{H^\frac{1}{4}}<\infty,
\end{equation*}
where the second inequality follows from \cref{l:43}. By Sobolev embedding,
\begin{equation*}
    \left|\left\{x:\varlimsup\limits_{t\rightarrow 0}\left|\int_0^t e^{i(t-t^\prime)\partial_{xx}}\Big(V_ke^{-it^\prime H}u_0\Big)dt^\prime\right|>0\right\}\right|=0,
\end{equation*}
for all $k\in \mathbb{Z}$ and this completes the proof.
\end{proof} 

\section{Quadratic Nonlinearities.}\label{qnonlinearity}
Consider the following qNLS Cauchy problem:
\begin{center}
    $\begin{cases}
    i\partial_t u + \partial_{xx}u = N_i(u,\overline{u}),\:u(0) = u_0 \in H^{s}(\mathbb{R}),\\
    N_1(u,\overline{u}) = u^2; N_2(u,\overline{u}) = u\overline{u}; N_3(u,\overline{u}) = \overline{u}^2.
    \end{cases}$
\end{center}

The well-posedness of qNLS above is studied in \cite{kenig1996quadratic}. By $X^{s,b}$ method, the qNLS for $N_1$ and $N_3$ are well-posed in $H^s(\mathbb{R})$ for $s>-\frac{3}{4}$ whereas that for $N_2$ is well-posedness for $s>-\frac{1}{4}$; the well-posedness associated to $N_1$ was improved to $H^{-1}(\mathbb{R})$ and was shown to be sharp in \cite{bejenaru2006sharp}. In the integral form, the solution satisfies
\begin{align*}
    u(t) = e^{it\partial_{xx}}u_0 -i \int_0^t e^{i(t-\tau)\partial_{xx}}N_i(u)(\tau)d\tau, t \in [-\delta,\delta],
\end{align*}
and the goal is to prove smoothing estimates for $N_i$, $i=1,2,3$ as in \cref{l:33} from which convergence to initial data follows by the Sobolev embedding.

\begin{lemma}\label{l:51}
Let $s \geq 0$, $a \in [0,\frac{1}{2})$ and $\delta\in (0,1]$. Then there exists $b=\frac{1}{2}+$ and $\gamma=\frac{1}{2}-$ such that $b<1-\gamma$ and the following estimates hold for $i=1,3$:
\begin{align*}
     \lVert N_i(u,\overline{u}) \rVert_{X^{s+a,-\gamma}_\delta} &\lesssim_{s,a,b,\gamma} \lVert u \rVert_{X^{s,b}_\delta}^2. 
\end{align*}
\end{lemma}

\begin{lemma}\label{l:52}
Let $s >\frac{1}{4}$, $a \in [0,\frac{1}{2}]$ and $\delta\in (0,1]$. Then there exists $b=\frac{1}{2}+$, $\gamma=\frac{1}{2}-$ such that $b<1-\gamma$ and the following estimate holds:
\begin{align*}
     \lVert N_2(u,\overline{u}) \rVert_{X^{s+a,-\gamma}_\delta} &\lesssim_{s,a,b,\gamma} \lVert u \rVert_{X^{s,b}_\delta}^2.   
\end{align*}
\end{lemma}

\begin{proof}[proof of \cref{l:51}]
The $N_3$-estimate will be shown to be an easy consequence of the $N_1$-estimate, and therefore we focus on the former. Denote
\begin{align*}
F(\xi,\tau) = |\hat{u}(\xi,\tau)|\langle \xi \rangle^s \langle \tau + \xi^2 \rangle^b;\: W(\xi,\tau,\xi_1,\tau_1) = \frac{\langle \xi \rangle^{s+a}\langle \tau + \xi^2 \rangle^{-\gamma}}{\langle \xi-\xi_1\rangle^s \langle \tau - \tau_1 + (\xi - \xi_1)^2\rangle^b\langle \xi_1 \rangle^s \langle \tau_1 + \xi_1^2 \rangle^b}.
\end{align*}
Neglecting $\delta$-dependence as before, we have
\begin{align*}
    \lVert u^2 \rVert_{X^{s+a,-\gamma}}^2 &= \left|\left| \int \langle \xi \rangle^{s+a} \langle \tau + \xi^2 \rangle^{-\gamma} \hat{u}(\xi-\xi_1,\tau-\tau_1)\hat{u}(\xi_1,\tau_1)d\xi_1 d\tau_1\right|\right|_{L^2_{\xi,\tau}}^2\\
    &\leq \left|\left| \int \frac{\langle \xi \rangle^{s+a}\langle \tau + \xi^2 \rangle^{-\gamma}}{\langle \xi-\xi_1\rangle^s \langle \tau - \tau_1 + (\xi - \xi_1)^2\rangle^b\langle \xi_1 \rangle^s \langle \tau_1 + \xi_1^2 \rangle^b} F(\xi-\xi_1,\tau-\tau_1)F(\xi_1,\tau_1)d\xi_1 d\tau_1\right|\right|_{L^2_{\xi,\tau}}^2\\
    &= \left|\left| \int W^2 d\xi_1d\tau_1 \cdot \int F(\xi-\xi_1,\tau-\tau_1)^2F(\xi_1,\tau_1)^2d\xi_1 d\tau_1\right|\right|_{L^1_{\xi,\tau}}\\
    &\leq \left|\left| \int W^2 d\xi_1 d\tau_1\right|\right|_{L^\infty_{\xi,\tau}} \cdot\lVert F^2 \ast F^2\rVert_{L^1_{\xi,\tau}}=\left|\left| \int W^2 d\xi_1 d\tau_1\right|\right|_{L^\infty_{\xi,\tau}}\cdot \lVert u \rVert_{X^{s,b}}^4.
\end{align*}

Hence, it suffices to prove that $\lVert \int W^2 d\xi_1 d\tau_1\rVert_{L^\infty_{\xi,\tau}}<\infty$.

By \cref{l:32}, we have
\begin{align*}
    \int \frac{d\tau_1}{\langle \tau - \tau_1 + (\xi - \xi_1)^2\rangle^{2b} \langle \tau_1 + \xi_1^2 \rangle^{2b}}\lesssim \langle \tau + (\xi_1-\xi)^2+\xi_1^2 \rangle^{-2b}.
\end{align*}

Similarly,
\begin{align*}
   \langle \tau + (\xi_1-\xi)^2+\xi_1^2 \rangle^{-2b}\langle \tau + \xi^2 \rangle^{-2\gamma}\lesssim \langle \xi_1(\xi_1-\xi)\rangle^{-2\gamma},
\end{align*}
and
\begin{align*}
    \frac{\langle \xi \rangle^{2s+2a}}{\langle \xi_1-\xi\rangle^{2s}\langle\xi_1\rangle^{2s}}\lesssim \langle \xi \rangle^{2a}.
\end{align*}

Altogether we have
\begin{align*}
    \sup\limits_{\xi,\tau}\int W^2 d\xi_1 d\tau_1\lesssim \sup\limits_{\xi}\Big(\langle \xi \rangle^{2a}\int \frac{d\xi_1}{\langle \xi_1(\xi_1-\xi)\rangle^{2\gamma}}\Big).
\end{align*}

Note that the integral is symmetric with respect to $\xi_1 = \frac{\xi}{2}$, and therefore $\int \frac{d\xi_1}{\langle \xi_1 (\xi_1-\xi)\rangle^{2\gamma}} = 2\int_{\xi/2}^{\infty} \frac{d\xi_1}{\langle \xi_1 (\xi_1-\xi)\rangle^{2\gamma}}$. Henceforth, assume $\xi \geq 0$ without loss of generality. On the region of integration, change variable $\eta = \xi_1 (\xi_1-\xi)= \xi_1^2 - \xi \xi_1$ to obtain:
\begin{align*}
     \xi_1 = \frac{\xi + \sqrt{\xi^2 + 4 \eta}}{2};\: d\xi_1 = \frac{d\eta}{\sqrt{\xi^2 + 4 \eta}}.
\end{align*}
and so the integral becomes
\begin{align*}
2\int_{\xi/2}^{\infty} \frac{d\xi_1}{\langle \xi_1 (\xi_1-\xi)\rangle^{2\gamma}}&= 2 \int_{-\xi^2/4}^\infty \frac{d\eta}{\sqrt{\xi^2 + 4 \eta}\langle \eta \rangle^{2\gamma}}= \int_0^\infty \frac{d\eta}{\sqrt{\eta} \langle \eta - \frac{\xi^2}{4}\rangle^{2\gamma}}. 
\end{align*}

Since this integral is bounded for all $\xi \in [0,1)$, it suffices to assume $\xi\geq 1$ and show $\int_0^\infty \frac{d\eta}{\sqrt{\eta} \langle \eta - \frac{\xi^2}{4}\rangle^{2\gamma}} \lesssim \frac{1}{\langle \xi \rangle^{4\gamma-1}}$. Then with $a <\frac{1}{2}$, it follows immediately that $\lVert \int W^2 d\xi_1 d\tau_1\rVert_{L^\infty_{\xi,\tau}}<\infty$, provided $b>\frac{1}{2}$ is chosen sufficiently small.\\

Let $c = \frac{\xi^2}{4}$ and estimate the integral in three different regions: i) $\eta \in [2c,\infty)$; ii) $\eta \in [\frac{c}{2},2c)$; iii) $\eta \in (0,\frac{c}{2})$.\\
\begin{align*}
    i) &: \int_{\eta \geq 2c} \lesssim \int_{2c}^\infty \frac{d\eta}{\sqrt{\eta} \eta^{2\gamma}} \simeq \frac{1}{\xi^{4\gamma-1}}.\\
    ii) &: \int_{\frac{c}{2}}^{2c} \lesssim \int_{\frac{c}{2}}^{2c} \frac{d\eta}{\sqrt{\eta} |\eta-c|^{2\gamma}} \lesssim c^{-1/2} \int_{\frac{c}{2}}^{2c} \frac{d\eta}{|\eta-c|^{2\gamma}} \lesssim c^{-1/2}\cdot c^{1-2\gamma} \simeq \frac{1}{\xi^{4\gamma-1}}.\\
    iii) &: \int_0^{\frac{c}{2}} \frac{d\eta}{\sqrt{\eta}\langle \eta -c \rangle^{2\gamma}} \lesssim \frac{1}{\langle c \rangle^{2\gamma}} \int_0^{\frac{c}{2}} \frac{d\eta}{\sqrt{\eta}} \lesssim \frac{c^{1/2}}{\langle c \rangle^{2\gamma}} \simeq \frac{1}{\xi^{4\gamma-1}}.
\end{align*}

Bringing all three cases together, we obtain the desired estimate, and this proves the first smoothing estimate.

As for the second estimate, for a general spacetime function $u$,
\begin{align*}
    \lVert \overline{u}\rVert_{X^{s.b}} = \lVert \hat{u}(\xi,\tau)\langle \xi \rangle^s \langle \tau - \xi^2 \rangle^b \rVert_{L^2_{\xi,\tau}}.
\end{align*}

Arguing as before, one obtains
\begin{align*}
    \lVert \overline{u}^2 \rVert_{X^{s+a,-\gamma}} \leq \left|\left| \int \Omega^2 d\xi_1 d\tau_1 \right|\right|_{L^\infty_{\xi,\tau}} \lVert u \rVert_{X^{s,b}}^4.
\end{align*}
where
\begin{align*}
    \Omega(\xi,\tau,\xi_1,\tau_1) = \frac{\langle \xi \rangle^{s+a}\langle \tau - \xi^2 \rangle^{-\gamma}}{\langle \xi-\xi_1\rangle^s \langle \tau - \tau_1 + (\xi - \xi_1)^2\rangle^b\langle \xi_1 \rangle^s \langle \tau_1 + \xi_1^2 \rangle^b},
\end{align*}
and therefore it suffices to show $\lVert \int \Omega^2 d\xi_1 d\tau_1 \rVert_{L^\infty_{\xi,\tau}}<\infty$. As before,\\
\begin{align*}
    \left|\left| \int \Omega^2 d\xi_1 d\tau_1 \right|\right|_{L^\infty_{\xi,\tau}}&\lesssim \left|\left| \int \frac{\langle \xi \rangle^{2s+2a}\langle \tau - \xi^2 \rangle^{-2\gamma}}{\langle \xi_1 - \xi \rangle^{2s}\langle \xi_1 \rangle^{2s} \langle \tau + (\xi_1-\xi)^2 + \xi_1^2 \rangle^{2b}} d\xi_1 \right|\right|_{L^\infty_\xi}\\
    &\leq \left|\left| \int \frac{\langle \xi \rangle^{2s+2a}}{\langle \xi_1-\xi \rangle^{2s}\langle \xi_1 \rangle^{2s}\langle \xi_1^2 -\xi_1\xi +\xi^2 \rangle^{2\gamma}}d\xi_1\right|\right|_{L^\infty_\xi}\lesssim \left|\left| \langle \xi \rangle^{2a} \int \frac{d\xi_1}{\langle \xi_1^2-\xi_1\xi +\xi^2\rangle^{2\gamma}}\right|\right|_{L^\infty_\xi},
\end{align*}
where these inequalities are direct applications of \cref{l:32}. Then by a direct computation,
\begin{align*}
    |\xi_1^2-\xi_1\xi +\xi^2|&\geq |\xi_1 (\xi_1-\xi)|.\\
    \Rightarrow\langle \xi_1^2-\xi_1\xi +\xi^2\rangle &\geq \langle \xi_1 (\xi_1-\xi)\rangle, \forall \xi,\xi_1 \in \mathbb{R}.
\end{align*}
Then $\lVert \int \Omega^2 d\xi_1 d\tau_1 \rVert_{L^\infty_{\xi,\tau}}<\infty$ follows from our previous result:
 \begin{align*}
     \left|\left| \langle \xi \rangle^{2a}\int \frac{d\xi_1}{\langle \xi_1 (\xi_1-\xi)\rangle^{2\gamma}}\right|\right|_{L^\infty_\xi}<\infty.
 \end{align*}
\end{proof}

\begin{proof}[proof of \cref{l:52}]
Arguing as before, it suffices to prove
\begin{align*}
\sup\limits_{\xi,\tau}\Big(\langle \xi \rangle^{2s+2a}\langle\tau+\xi^2\rangle^{-2\gamma}\int \frac{d\xi_1d\tau_1}{\langle \xi-\xi_1\rangle^{2s}\langle \xi_1\rangle^{2s}\langle \tau_1-(\tau+(\xi-\xi_1)^2)\rangle^{2b}\langle \tau_1-\xi_1^2\rangle^{2b}}\Big)<\infty.    
\end{align*}
For $|\xi|<1$,
\begin{align*}
    \langle \xi \rangle^{2s+2a}\langle\tau+\xi^2\rangle^{-2\gamma}\int \frac{d\xi_1d\tau_1}{\langle \xi-\xi_1\rangle^{2s}\langle \xi_1\rangle^{2s}\langle \tau_1-(\tau+(\xi-\xi_1)^2)\rangle^{2b}\langle \tau_1-\xi_1^2\rangle^{2b}}&\lesssim \int \frac{d\xi_1}{\langle \xi_1 \rangle^{4s}\langle \tau + \xi^2 - 2 \xi \xi_1\rangle^{2b}}\\
    &\lesssim_{s,a}\int \frac{d\xi_1}{\langle \xi_1 \rangle^{4s}}\leq C<\infty,
\end{align*}
where the upper bound $C$ is independent of $\tau$. For $|\xi|\geq 1$, changing variable $z=2\xi\xi_1-(\tau+\xi^2)$,
\begin{align*}
    \langle \xi \rangle^{2s+2a}\langle\tau+\xi^2\rangle^{-2\gamma}\int \frac{d\xi_1d\tau_1}{\langle \xi-\xi_1\rangle^{2s}\langle \xi_1\rangle^{2s}\langle \tau_1-(\tau+(\xi-\xi_1)^2)\rangle^{2b}\langle \tau_1-\xi_1^2\rangle^{2b}}&\lesssim \frac{\langle \xi \rangle^{2a}}{\langle \tau+\xi^2\rangle^{2\gamma}}\int \frac{d\xi_1}{\langle 2\xi\xi_1-(\tau+\xi^2)\rangle^{2b}}\\
    &\simeq \frac{\langle \xi \rangle^{2a-1}}{\langle \tau+\xi^2\rangle^{2\gamma}}\int\frac{dz}{\langle z \rangle^{2b}}\lesssim \langle \xi \rangle^{2a-1}.
\end{align*}
\end{proof}
\begin{remark}
As for the smoothing estimate for $N_2$, the condition $s>\frac{1}{4}$ is necessary to make certain integrals converge; in fact if $\xi=\tau=0$, then the expression inside the $\sup\limits_{\xi,\tau}$ (see the proof for \cref{l:52}) is
\begin{align*}
    \int \frac{d\xi_1 d\tau_1}{\langle \xi_1 \rangle^{4(\frac{1}{4})}\langle \tau_1-\xi_1^2\rangle^{4b}}=\infty.
\end{align*}
\end{remark}

\begin{proof}[proof of \cref{main2}]
The positive statements are consequences of \cite{carleson1980some} and Duhamel nonlinear terms being continuous in space and time via the smoothing estimates followed by the Sobolev embedding. We focus on the negative part.

For $s\in (0,\frac{1}{4})$ we know from \cite{dahlberg1982note} that there exists $u_0\in H^s(\mathbb{R})$ such that convergence to initial data fails on some set $E$ of positive measure. By Lemma $3.1$, we choose $a=\frac{1}{2}-$ to obtain
\begin{align*}
    DN(x,t)\coloneqq\int_0^t e^{i(t-t^\prime)\partial_{xx}}N_i(u,\overline{u})(t^\prime)dt^\prime \in C^0_t H^{\frac{1}{2}+}_x ([0,\delta]\times\mathbb{R}).
\end{align*}
By the triangle inequality,
\begin{align*}
    |u(t)-u_0|\geq |e^{it\partial_{xx}}u_0-u_0|-|DN(x,t)|.
\end{align*}

By continuity, $\lim\limits_{t\rightarrow 0}|DN|=0$ a.e., and therefore
\begin{align*}
    \left|\left\{x\in E: \varlimsup\limits_{t\rightarrow 0}|u(t)-u_0|>0\right\}\right|\geq \left|\left\{x\in E: \varlimsup\limits_{t\rightarrow 0}|e^{it\partial_{xx}}u_0-u_0|>0\right\}\right|>0. 
\end{align*}

Since $H^s(\mathbb{R})\hookrightarrow L^2(\mathbb{R})$ for $s \in (0,\frac{1}{4})$, a.e. pointwise convergence cannot hold for initial data in $L^2(\mathbb{R})$, and this finishes the proof.
\end{proof}

\section{Negative Result: Baire Category Approach.}\label{negative2}
Note that if a.e. pointwise convergence does not hold for $s<\frac{1}{4}$, then it also fails for $0\leq s^\prime \leq s$. Define $D(s)$ to be the collection of $f\in H^s(\mathbb{R})$ with a compact support such that $\varlimsup\limits_{t\rightarrow 0}|e^{-itH}f|\gtrsim 1$ uniformly on some measurable set (of positive measure) $E_f\subseteq\mathbb{R}\setminus supp(f)$. Define $D_0(s)$ similarly via $e^{it\partial_{xx}}$. One motivation for considering functions in $D(s)$ comes from Sj\"olin's work on localisation of Schr\"odinger means.  
\begin{lemma}\cite{sjolin2013nonlocalization} Let $s<\frac{1}{4}$. There exists $f = \sum\limits_{n=1}^\infty f_n \in H^s(\mathbb{R})$ supported in $(-\frac{\delta}{4},\frac{\delta}{4})$ for some $\delta>0$ where $f_n$'s are smooth and $\varlimsup\limits_{t\rightarrow 0} |e^{it\partial_{xx}}f(x)|\geq c>0$ uniformly on a measurable set $E \subseteq (\frac{\delta}{2},\delta)$ of positive measure.
\end{lemma}
\begin{remark}
Given $K\subseteq \mathbb{R}$, a compact subset, one can modify the arguments of \cite{sjolin2013nonlocalization} to explicitly construct $f \in H^s(\mathbb{R})$ with its support in $K$ such that $e^{it\partial_{xx}}f \rightarrow f$ as $t\rightarrow 0$ fails in a.e. sense on $K^c$.
\end{remark} 
We show that $D_0(s) = D(s)$ for $s<\frac{1}{4}$, or i.e., that in the short-time limit, the potentials play no role in the convergence of solutions.
\begin{proposition}\label{p:62}
Let $V\in L^2(\mathbb{R})$ and $0<s<\frac{1}{4}$. Then, $D_0(s) = D(s)$. 
\end{proposition}
\begin{proof}[proof of \cref{p:62}]
Writing $u(t)=e^{-itH}f$, the Duhamel formula yields
\begin{equation*}
    u(t)-f = e^{it\partial_{xx}}f-f -i\int_0^t e^{i(t-\tau)\partial_{xx}}(Vu)(\tau)d\tau.
\end{equation*}
As in the proof of Theorem $2.1$, we apply the smoothing estimate (\cref{l:33}) on $Vu$ by choosing $a=\frac{1}{2}-$ and the well-posedness result (\cref{l:31}) to obtain that the Duhamel integral term is continuous in time and $H^{\frac{1}{2}+}$ in space, from which $D_0(s)=D(s)$ follows immediately.
\end{proof} 
\begin{remark}
In the proof, note that our smoothing estimate is insufficient to conclude $D_0(0)=D(0)$.
\end{remark}
Motivated by \cref{p:62}, we restrict the collection of counterexamples to $f\in D(s)=D_0(s)$, assuming $V \in L^2$, and therefore $e^{-itH}$ can be replaced by $e^{it\partial_{xx}}$. Fix $J=(-1,1), \phi \in C^\infty_c(K)$ where $K\subseteq \mathbb{R}\setminus J$ is compact. It turns out that it is not an easy task to explicitly find such examples. Another more commonly-used approach is via the Stein-Nikisin maximal principle (\cite{nikishin1972resonance}), which states:
\begin{lemma}
$e^{it\Delta}f\xrightarrow{a.e.}f$ as $t\rightarrow 0$ for all $f\in H^s(\mathbb{R}^n)$ if and only if
\begin{equation}\label{Nikisin1}
\left|\left| \sup_{0<t<1} |e^{it\Delta}f|\right|\right|_{L^2(B(0,1))}\lesssim_{s,n} \lVert f \rVert_{H^s(\mathbb{R}^n)}, \forall f \in C^\infty_c(\mathbb{R}^n).    
\end{equation}
\end{lemma}
Note that the $L^2(B(0,1))$ cannot be strengthened to $L^p(B(0,1))$ for $p\in (2,\infty]$ when $s<\frac{1}{4},\:n=1$ due to the H\"older's inequality. For such $p$, we ask whether the $f$ on the LHS of \cref{Nikisin1} can be replaced by $f\phi$, i.e., whether the $H^s$ norm controls the $\phi$-localised maximal operator. It turns out that this fails for a big class of functions.
\begin{proposition}\label{p:71}
For $s<\frac{1}{4}$ and $p\in (2,\infty)$, the following strong-type estimate fails:
\begin{equation}\label{Nikisin2}
    \left|\left| \sup_{0<t<1}|e^{it\partial_{xx}}(f\phi)|\right|\right|_{L^p(J)}\lesssim \lVert f \rVert_{H^s(\mathbb{R})}.
\end{equation}
\end{proposition} 
Note that if \cref{Nikisin2} fails for $p$, then it fails for $\Tilde{p}\geq p$. On the other hand, Sj\"olin in \cite{sjolin2012some} showed that for every $f\in H^s(\mathbb{R}^n)$ with a compact support, $e^{it\Delta}f\xrightarrow[t\rightarrow 0]{} 0$ \textit{for all} $x\in \mathbb{R}^n\setminus supp(f)$ if and only if $s\geq \frac{n}{2}$. Since the free Schr\"odinger operator is given by the convolution $e^{it\Delta}f = K_t \ast f$ where $K_t(x) = (4\pi it)^{-\frac{n}{2}}e^{i\frac{|x|^2}{4t}}$, it is evident that $e^{it\Delta}f\in C^\infty_x(\mathbb{R}^n)$ for each $t\in \mathbb{R}\setminus \left\{0\right\}$ since $K_t \in C^\infty_x(\mathbb{R})$ and $f$ has a compact support, and hence it makes sense to evaluate $e^{it\Delta}f$ pointwise. Sj\"olin showed, via Baire category approach, that for $s<\frac{n}{2}$ there exists $f\in H^s(\mathbb{R}^n)$ with a compact support in $S=\left\{|x|\in (1,2)\right\}$ such that $e^{it\Delta}f(0)\rightarrow \infty$ as $t\rightarrow 0$. Hence $\lVert e^{it\Delta}f\rVert_{L^\infty (B(0,1))}\rightarrow \infty$ as $t\rightarrow 0$ since $e^{it\Delta}f$ is smooth. Here we are interested in the $L^p$-behavior of solutions in the short-time limit. For $p\in [1,2]$, $\lVert e^{it\Delta}f\rVert_{L^p(B(0,1))}$ stays bounded due to the $L^2$-conservation of solutions. For $p\in (2,\infty)$, it is unclear whether the solution blows up or stays bounded. We show a weaker result that the $L^p$-norm of solutions diverges in some time-averaged sense:
\begin{proposition}\label{p:72}
Let $\left\{t_k\right\}_{k=1}^\infty$ be a real-sequence contained in $(0,1]$ that tends to zero as $k\rightarrow \infty$ and $p\in (2,\infty)$. Then there exists a dense, $G_\delta$ residual set $\mathcal{C}\subseteq H^s(\mathbb{R})$ such that for every $f\in\mathcal{C}$, $\left\{e^{it_k\partial_{xx}}(f\phi)\right\}\notin l^q L^p(\mathbb{N}\times J)$ for all $q\in [1,p]$.\footnote{Recall that a measurable set is $G_\delta$ if it can be realised as a countable intersection of open sets. A set is \textit{meager} if it can be realised as a countable union of nowhere dense sets, and its complement is called a \textit{residual}.}
\end{proposition}
Our proof is a simple application of the Banach-Steinhaus theorem. Given a sequence $ t_k\xrightarrow[k\rightarrow \infty]{}0$, define 
\begin{equation*}
S_n f = \sup\limits_{k\leq n}|e^{it_k\partial_{xx}}(f\phi)|,\: S f =\sup\limits_{k}|e^{it_k\partial_{xx}}(f\phi)|,    
\end{equation*}
for $f\in \mathscr{S}(\mathbb{R})$. It is straightforward to verify
 \begin{equation}\label{sublinearity}
     S_n(f+g)\leq S_nf+S_n g,\: S_n(\lambda f) = |\lambda| S_n f,\forall \lambda \in \mathbb{C},
 \end{equation}
and in particular, $S_n$ is not linear. For this reason, the following extension of the Banach-Steinhaus Theorem, traditionally studied in the context of linear operators, is applied where its proof could be done as \cite[Theorem 5.8]{rudin2006real}.
\begin{lemma}\label{l:72}
Let $\left\{T_\alpha\right\}_{\alpha\in A}$ be a family of continuous operators on $X$ into $Z=L^p(Y,\nu)$ for $p\in [1,\infty]$ where $X$ is a Banach space, $(Y,\nu)$ is a $\sigma$-finite measure space and $A$ is some directed set.
\begin{equation}\label{Sublinear}
    \lVert T_\alpha (x+y)\rVert \leq \lVert T_\alpha x \rVert + \lVert T_\alpha y \rVert,\: \lVert T_\alpha (\lambda x)\rVert = \lambda \lVert T_\alpha x \rVert,\forall \lambda\geq 0.
\end{equation}
Then either $\lim \limits_{x \rightarrow 0} \lVert T_\alpha x \rVert = 0$ uniformly in $\alpha$, i.e., $\left\{T_\alpha\right\}$ is equicontinuous at the origin, or\\ $\left\{ x \in X: T_\alpha x\; \text{is unbounded in Z} \right\}$ forms a residual set that is dense $G_\delta$ in $X$.
\end{lemma} 

\begin{proof}[proof of \cref{p:71,p:72}]
We first claim that $\left\{S_n\right\}$ defines a family of continuous sublinear operators on $H^s(\mathbb{R})$ into $L^p(J)$ that satisfies the hypotheses of \cref{l:72}. By the triangle inequality,
\begin{equation*}
    |S_n f- S_ng|\leq S_n(f-g).
\end{equation*}
Hence it suffices to show that $S_n$ is a bounded map to show continuity. Since $l^p(\mathbb{N})\hookrightarrow l^\infty(\mathbb{N})$, we obtain
\begin{align}\label{Sublinear2}
    \lVert S_n f\rVert_{L^p(J)} &= \lVert e^{it_k\partial_{xx}}(f\phi)\rVert_{L^p(J) l^\infty_{k\leq n}}\leq \lVert e^{it_k\partial_{xx}}(f\phi)\rVert_{ l^p_{k\leq n}L^p(J)}\\
    &\lesssim \left|\left| |t_k|^{-(\frac{1}{2}-\frac{1}{p})}\lVert f \phi\rVert_{L^{p^\prime}(J)}\right|\right|_{l^p_{k\leq n}} \lesssim_{p,\phi,n} \lVert f \rVert_{H^s(\mathbb{R})},\nonumber
\end{align}
and hence the continuity. From \cref{sublinearity}, \cref{Sublinear} could be verified. It is shown, by contradiction, that the $\left\{S_n\right\}$ cannot be equicontinuous at the origin. Assume it is. Then $S$ is continuous in measure at the origin. Suppose $f_j\rightarrow 0$ in $H^s(\mathbb{R})$ as $j\rightarrow \infty$ and let $\lambda>0$. Let $\epsilon>0$ for which there exists $\delta>0$ such that $\lVert S_n f\rVert_{L^p(J)}<\epsilon^{\frac{1}{p}}\lambda$ for all but finitely many $n\in\mathbb{N}$ and all $f$ such that $\lVert f \rVert_{H^s}<\delta$. Then let $j \geq N$, some $N\in \mathbb{N}$ sufficiently big, such that $\lVert f_j\rVert_{H^s}<\delta$ for all $j\geq N$, and let $n$ be sufficiently large such that $\left|\left\{|Sf_j -S_nf_j|\geq  \frac{\lambda}{2}\right\}\right|\lesssim \epsilon$; recall that $S_n f_j \xrightarrow{a.e.}Sf_j$ implies $S_nf_j \rightarrow Sf_j$ in measure on a finite measure space as $n\rightarrow\infty$. Then we obtain
\begin{equation*}
    \left|\left\{|S f_j|\geq \lambda\right\}\right|\leq \left|\left\{|Sf_j -S_n f_j|\geq \frac{\lambda}{2}\right\}\right|+\left|\left\{|S_n f_j|\geq \frac{\lambda}{2}\right\}\right|\lesssim \epsilon,
\end{equation*}
where the second term is bounded above by $\epsilon$ up to a constant by the Chebyshev's inequality.

Now it is shown that convergence a.e. to initial data holds for all $f\in H^s(\mathbb{R})$ with a compact support, which is a contradiction since $s<\frac{1}{4}$ and due to the explicit construction of an initial datum with a compact support in \cite{sjolin2013nonlocalization}. Pick $f_n\rightarrow f\in H^s(\mathbb{R})$ where $f_n \in \mathscr{S}(\mathbb{R})$. Then,
\begin{equation*}
\begin{split}
\left|\left\{x\in J:\varlimsup\limits_{k\rightarrow\infty}|e^{it_k\partial_{xx}}(f\phi)|>\lambda\right\}\right|&=\left|\left\{x\in J:\varlimsup\limits_{k\rightarrow\infty}|e^{it_k\partial_{xx}}((f-f_n)\phi)|>\lambda\right\}\right|\\
    &\leq\left| \left\{x\in J:S (f-f_n)>\lambda\right\}\right|\xrightarrow[n\rightarrow 0]{} 0,    
\end{split}
\end{equation*}
for all $\lambda>0$ where the last limit follows from the continuity in measure of $S$. Hence the supposed equicontinuity fails and there exists a dense $G_\delta$ set $\mathcal{C}\subseteq H^s(\mathbb{R})$ such that if $f \in \mathcal{C}$, then $\left\{S_n f \right\}$ is unbounded in $L^p(J)$. By monotonicity, $\lVert S_n f\rVert_{L^p(J)}\leq \lVert Sf\rVert_{L^p(J)}\leq \left|\left| \sup\limits_{0<t<1}|e^{it\partial_{xx}}(f\phi)|\right|\right|_{L^p(J)}$, and therefore, \cref{Nikisin2} cannot hold for every $f\in \mathcal{C}$. By the right-most estimate in \cref{Sublinear2}, we obtain that $\left\{e^{it_k\partial_{xx}}(f\phi)\right\}\notin l^p L^p(\mathbb{N}\times J)$ for all $f\in \mathcal{C}$.
\end{proof}

\section{Acknowledgements.}

The author would like to appreciate his doctoral advisor Mark Kon for insightful comments.

\bibliographystyle{abbrv}
\bibliography{ref2}
\end{document}